\documentstyle[12pt]{article}
\title{Kirby-Melvin's  $\tau_r^{'}$ and Ohtsuki's $\tau$ for
 Lens Spaces
 }
\author{ Bang-He Li \& Tian-Jun Li\thanks{The first author is supported partially by the Tianyuan Foundation of
P. R. China, the second author is supported partially by NSF grant DMS
9304580.}}
\date{}
\begin{document}
\maketitle
\baselineskip 18pt
\begin{abstract}

Explicit formulae for $\tau_r^{'}(L(p,q))$  and $\tau(L(p,q))$  are obtained
for all $L(p,q)$.
 \end{abstract}
\vskip 36pt

There are three systems of invariants of Witten-type for closed
oriented 3-manifolds:
\bigskip
\begin{description}
\item \quad 1.~$\{\tau_r(M), r\geq2; \tau_r^{'}(M),{\hskip .1cm} r{\hskip .1cm}
\mbox{odd}\geq3\}$, \mbox{where} $\tau_r$ was defined by Reshetikhin
and Turaev [1], and $\tau_r^{'}$ \mbox{was defined by Kirby-Melvin [2]}
\smallskip
\item \quad 2.~$\{\Theta_r(M,A), r\geq1, \mbox{where} A \,\,\mbox{is a 2r-th
primitive root
of unity}\}$ {defined by Blanchet, Habegger, Masbamm and Vogel  [3]}
\smallskip
\item \quad 3.~$\{\xi_r(M,A), r\neq1,4,4k+2,\quad \mbox{where} A \,\,\mbox{is
an r-th primitive root of unity}\}$ defined by the first author [4].
\end{description}
\smallskip
And it was proved in [4] that they are equivalent.
Explicit formulae for ~$\{\tau_r(M), r\geq2\}$ and ~$\{\xi_r(M,e_r),
{\hskip .1cm} r{\hskip .1cm}\mbox{odd}\geq3\}$ have been obtained for
lens spaces in [5], where $e_a=\exp(2\pi\sqrt{-1}/ a)$.

Ohtsuki [6] defined his invariant
$$\tau(M)=\sum_{n=0}^{\infty}\lambda_n(t-1)^n \in Q[[t-1]]$$
for rational homology 3-sphere $M$, and obtained that
$$\tau(L(p,q))=t^{-3s(q,p)}{{t^{1\over 2p}-t^{-{1\over 2p}}}\over
{t^{1\over 2}-t^{-{1\over 2}}}}$$
for lens spaces $L(p,q)$ with $p$ odd, where $s(q,p)$ is the Dedekind sum.

To obtain $\tau(L(p,q))$ with $p$ odd, Ohtsuki used the formula for
$\tau^{'}_r(L(p,q))$ with $r$ an odd prime, $p$ odd and not divisible by
$r$, found by Kirby-Melvin [7], and Garoufalidis [8].

To obtain\ $\tau(L(p,q))$ with $p$ even, explicit formulae for ~$\tau_r^{'}$
in this case are needed.

The aim of this note is to first derive explicit formulas of $\tau_r^{'}$
for all $L(p,q)$ and all $r$, and then give formulas of $\tau(L(p,q))$.

Notice that in [6], the orientation of $L(p,q)$ comes from the continued fraction
expantion of $p/q$, while in [5], the orientation comes from that
 of $-p/q$. So to calculate $\tau^{'}_r(L(p,q))$  with orientation as
 in [6], we should start from the conjugate of the formula in Theorem 4.1 in
 [5], that is $\overline{\xi_r(L(p,q))}=\xi_r(L(p,q), e_r)$.

Let $a,b$ and $r$ be integers, denote by $(b,r)$ the greatest common divisor of
$b$ and $r$. If $(b,r)=1$, denote $(a/b)^{\surd}=ab' \in Z/{rZ}$, where $b'b
\equiv 1 \pmod r$. Notice that if $a/b=a_1/b_1$ and $(b_1,r)=1$, then
$ab'=a_1b_1^{'}$ in $Z/{rZ}$.

{\bf Theorem.} Let $r>1$ be odd and $C=(p,r)$, then
$$\tau_r^{'}(L(p,q))=\cases{\displaystyle({p\over r})e_{r}^{-(3s(q,p))^{\surd}}
{e_r^{2p^{'}}-e_r^{-2p^{'}}\over e_r^{2^{'}}-e_r^{-2^{'}}},&if $c=1$ \cr
\quad\cr
\quad\cr
\displaystyle (-1)^{{r-1\over2}{c-1\over2}}({{p/ c}\over{r/ c}})
({q(1\mp r)/4\over c})
(e_{r}^{-12s(q,p)}
\cr
\displaystyle e_{pc}^{-({r/c})^{'}(q+q^*-\eta p^*p)}
e_{rc}^{2\eta({p/c})^{'}})^{{1\mp r}\over 4}
{{\epsilon(c)\sqrt c\eta}\over {e_r^{-2^{'}}-e_r^{2^{'}}}} ,
&if  $c>1$ ,  $c\mid q^*+\eta$  \cr
\quad& and $r\equiv \pm 1 \pmod 4$\cr
\quad\cr
0,&if $c>1$ and $c\mid\!\llap /  q^{*}\pm 1$\cr}$$
where $\eta=1$ or $-1$, $p^*p+q^*q=1$ with $0< q^{*} <p$,
$2'2\equiv 2p(2p)'\equiv 1 \pmod r$, $({a\over b})$ is the Jacobi symbol,
 $({p/ c})^{'}{p/ c}
+ ({r}/{c})'{{r}/{c}}= 1$, and
$$\epsilon (c)=\cases {1& if $c\equiv 1 \pmod 4$\cr
i & if $c\equiv -1 \pmod 4$\cr}$$
\vspace{.5cm}

Proof. By the formulas in [4] (see also [5]),
we have

$$\tau_r^{'}(M)={(\sqrt{4\over r}sin{\pi\over r})}^{\nu}
\Theta_{r}(M, \pm ie_{4r})$$ for $ r\equiv \pm 1\pmod 4$
and
$$\xi_r(M,A)=
 2^{\nu}\Theta_r(M, -A)
$$ for $ 3\leq r\equiv 1\pmod2,$
where $\nu$ is the first Betti number. Since for lens spaces, $\nu=0$,
we have,  for $r\equiv \pm 1\pmod 4$,
\bigskip
$$\begin{array}{cl}
{\tau_r^{'}(L(p,q))}&=\displaystyle{\Theta_r(L(p,q), \pm ie_{4r})}\\
&=\xi_r (L(p,q), \mp ie_{4r})\\
&=\xi_r (L(p,q),e_r^{{1 \mp r}\over 4}),
\end{array}$$

{\bf Case 1.} $c=1$. By Theorem 4.1 in [5]
$$                         \xi_r (L(p,q), e_{r})
=\displaystyle{(\frac{p}{r})e_{r}^{-12s(q,p)}e_{p}
^{r'(q+q^{*})}
\frac{e_{r}^{2p'}-e_{r}^{-2p'}}{e_{r}^{2}-e_{r}^{-2}}}.$$

Since $e_{r}^{-12s(q,p)}e_{p}^{r'(q+q^{*})}=e^a_r$ for some $a \in Z$ (by
checking the proof of Theorem 4.1 in [5]), it follows that
$$\tau_r^{'}(L(p,q))=
\displaystyle (\frac{p}{r})e_{r}^{-12s(q,p){{1\mp r}\over 4}}e_{p}
^{r'(q+q^{*}){{1\mp r}\over 4}}
\frac{e_{r}^{{{1\mp r}\over 4}2p'}-e_{r}^{-{{1\mp r}\over 4}2p'}}{e_{r}^
{{{1\mp r}\over 4}2}-e_{r}^{-{{1\mp r}\over 4}2}}$$
for $r\equiv \pm 1 \pmod 4.$  
The simple relations ${{1\mp r}\over 2}=2'$ and $2'p'=(2p)'$ in
$Z/ {rZ}$ imply that
$$\frac{e_{r}^{{{1\mp r}\over 4}2p'}-e_{r}^{-{{1\mp r}\over 4}2p'}}{e_{r}^
{{{1\mp r}\over 4}2}-e_{r}^{-{{1\mp r}\over 4}2}}
=
\frac{e_{r}^{2p'}-e_{r}^{-2p'}}{e_{r}^{2'}-e_{r}^{-2'}}.$$

It is well known that (cf. [8] or [5])
$$-12s(q,p)= b- {{q+q^*}\over p}\; \;\mbox{for some}\; b\in Z.$$

Hence $3s(q,p)$ has a form of $m/{4p}$ for some $m\in Z$. Now $(4p,r)=1$, so
there are integers $P$ and $R$ such that
$ 4pP+rR=1.$
Thus $$\begin{array}{cl}{\displaystyle
e_{r}^{-12s(q,p){{1\mp r}\over 4}}e_{p}
^{r'(q+q^{*}){{1\mp r}\over 4}}}&
=
{\displaystyle
e_{r}^{-{m\over 4p}(4pP+rR)
{(1\mp r)}}
e_{p}
^{r'(q+q^{*}){{1\mp r}\over 4}}}\\
&=\displaystyle{e_{r}^{-{mP}}e_r^{
-{{1\mp r}\over 4}
(12s(q,p)rR-rr'{{q+q^{*}}\over p})}}.
\end{array}$$
Now
$$12s(q,p)rR-rr'{{q+q^{*}}\over p}=(-b+{{q+q^{*}}\over p})rR-
rr'{{q+q^{*}}\over p}=-brR+r(R-r'){{q+q^{*}}\over p}.$$
Since
$r'\equiv R\equiv 1 \pmod p$, $p| R-r'$. It is concluded that
$$12s(q,p)rR-rr'{{q+q^{*}}\over p}\equiv 0 \pmod r.$$

Because $mp=(3s(q,p))^{\surd} $, the proof is complete for case 1.

{\bf Case 2.} $c>1$ and $c| q^*+\eta$ for $\eta=1$ or $-1$. Then
$$\xi_r(L(p,q), e_r)=
\displaystyle (-1)^{{r-1\over2}{c-1\over2}}({{p/ c}\over{r/ c}})({q\over c})
e_{r}^{12s(q,p)}
\displaystyle e_{pc}^{-({r/ c})^{'}(q+q^*-\eta p^*p)}
e_{rc}^{2\eta({p/c})^{'}}
{\epsilon(c)\sqrt c\eta\over e_r^{-2}-e_r^{2}}.$$

Since $c$ is odd and $({1\over c})=1$, it follows from Theorem 2.2 in [5] that
$$\epsilon(c)\sqrt c=\sum_{j=1}^c e_c^{j^2}=\sum_{j=1}^c e_r^{{r\over c}j^2}.$$
Hence, by the well-known formula for Gaussian sum (cf. [10] or [5]), if
$r\equiv \pm 1\pmod 4$,
$$\sum_{j=1}^c e_r^{{{1\mp r}\over 4}{r\over c}j^2}=
\sum_{j=1}^c e_c^{{{1\mp r}\over 4}j^2}=\epsilon(c)(\frac{(1\mp r)/4}{c})
\sqrt c.
$$
By the proof of Theorem 4.1 in [5],
$$e_{r}^{12s(q,p)}
\displaystyle e_{pc}^{-({r/c})^{'}(q+q^*-\eta p^*p)}
e_{rc}^{2\eta({p/c})^{'}}=e_r^k$$
for some $ k\in Z$.
Again from
$$e_r^{{{1\mp r}\over 4}2}=e_r^{{1\mp r}\over 2}=e_r^{2'},$$
it is concluded that  if $r\equiv \pm 1 \pmod 4$
$$\tau_r^{'}(L(p,q))= \xi_r (L(p,q), e_{r}^{{{1\mp r}\over 4}})=
\displaystyle (-1)^{{r-1\over2}{c-1\over2}}({{p/ c}\over{r/ c}})
({q(1\mp r)/4\over c})
e_r^{{{1\mp r}\over 4}k}
{{\epsilon(c)\sqrt c\eta}\over {e_r^{-2'}-e_r^{2'}}}.$$

{\bf Case 3.} $c>1$ and $c\mid\!\llap /  q^{*}\pm 1$. Since
$    \xi_r (L(p,q), e_{r})=0$, $\tau_r^{'}(L(p,q))=0 $. The theorem is proved.

{\bf Corollary}. For any lens spaces $L(p,q)$, Ohtsuki's invariant
$$\tau(L(p,q))=t^{-3s(q,p)}{{t^{1\over 2p}-t^{-{1\over 2p}}}\over
{t^{1\over 2}-t^{-{1\over 2}}}}.$$

Proof. Prop.5.2 in [6] for $p$ odd is now also true for  $p$ even by the
above theorem. Hence for even $p$, $\tau (L(p,q))$ has the same form as for odd $p$.

{\bf Remark 1}. Lemma 5.4 in [6] concerning odd $p$ is also true for even $p$.
Therefore, Remark 5.3 in [6]  for odd $p$ is true all $p$.

 {\bf Remark 2}. In [6], Ohtsuki regarded $\tau^{'}_r$ as $SO(3)$ invariants,
 while in [11], $SO(3)$ invariants were referred to those drived from the
 Kauffman module $Z[A, A^{-1}][z^2]$. The latter includes the former and
 something more.

 \begin{center} {\bf\Large \bf   References}\end{center}
 \begin{description}
\item{} [1] Reshetikhin, N.Yu., Turaev, V.G.: { Invariants of 3-manifolds
via link polynomials and quantum groups}, Invent. Math. {\bf 103} (1991), 
547-597
\item{} [2]    Kirby, R., Melvin, P.: { The 3-manifold invariants of Witten
and Reshetikhin-Turaev for $sl(2,C)$}, Invent. Math. {\bf 105} (1991), 473-545
\item{} [3]   Blanchet, C., Habegger, N., Masbaum, G., Vogel, P.: { Three-
manifold invarinats derived from the Kauffman bracket}, Topology, {\bf 31}
(1992), 685-699
\item{} [4] Li, B.H., { Relations among Chern-Simons-Witten-Jones invariants
}, Science in China, series A, {\bf 38} (1995), 129-146
\item{} [5] Li, B.H., Li, T.J.:{ Generalized Gaussian Sums and Chern-Simons-
Witten-Jones invariants of Lens spaces}, J. Knot theory and its Ramifications,
vol.5 No. 2 (1996) 183-224
\item{} [6] Ohtsuki, T., { A polynomial invariant of rational homology
3-spheres}, Inv. Math. {\bf 123} (1996), 241-257
\item{} [7] Kirby, R., Melvin, P.: { Quantum invariants of Lens spaces and
a Dehn sugery formula}, Abstract Amer. Math. Soc. {\bf 12} (1991), 435
\item{} [8]      Garoufalidis, S., {Relation among 3-manifold invariants},
Ph.D. Thesis, Univ. of Chicago, 1992
\item{} [9] Hickerson, D.,{ Continued fraction and density results}, J. Reine.
Angew. Math. {\bf 290} (1977), 113-116
\item{} [10]    Lang, S., { Algebraic Number Theory}, Springer, New York,
1986
\item{} [11] Li, B.H., Li, T.J.:{ SO(3) three-manifold invariants from
the Kauffman bracket}, Proc. of the Conference on Quantum Topology, Kansas,
1994, 247-258
\end{description}
\vspace{.3cm}

 Author's address
$$\begin{tabular}{ll}
\mbox{ Bang-He Li                   }& \mbox{ Tian-Jun Li } \\
            \mbox{ Institute of Systems Science,  }&
    \mbox{ School of math. } \\
            \mbox{ Academia Sinica }&
    \mbox{ IAS } \\
    \mbox { Beijnig 100080}& \mbox{ Princeton NJ 08540 }  \\
    \mbox { P. R. China }& \mbox{ U.S.A.} \\
    \mbox{ Libh@iss06.iss.ac.cn} &
    \mbox{ Tjli@IAS.edu}\\
\end {tabular}$$
\end{document}